\documentclass[%
 aip,
 amsmath,amssymb,
 reprint,%
]{revtex4-1}

\usepackage{graphicx}
\usepackage{dcolumn}
\usepackage{bm}

\usepackage{float}

\usepackage[utf8]{inputenc}
\usepackage[T1]{fontenc}
\usepackage{mathptmx}
\usepackage{etoolbox}

\makeatletter
\def\@email#1#2{%
 \endgroup
 \patchcmd{\titleblock@produce}
  {\frontmatter@RRAPformat}
  {\frontmatter@RRAPformat{\produce@RRAP{*#1\href{mailto:#2}{#2}}}\frontmatter@RRAPformat}
  {}{}
}%
\makeatother
\begin{document}

\preprint{AIP/123-QED}

\title[]{Dynamics, interference effects and multistability in a Lorenz-like system of a classical wave-particle entity in a periodic potential}
\author{J. Perks}
\author{R. N. Valani}%
 \email{rahil.valani@adelaide.edu.au}
\affiliation{School of Computer and Mathematical Sciences, University of Adelaide, South Australia 5005, Australia
}%

\date{\today}

\begin{abstract}

A classical wave-particle entity (WPE) can be realized experimentally as a droplet walking on the free surface of a vertically vibrating liquid bath, with the droplet's horizontal walking motion guided by its self-generated wave field. These self-propelled WPEs have been shown to exhibit analogs of several quantum and optical phenomena. Using an idealized theoretical model that takes the form of a Lorenz-like system, we theoretically and numerically explore the dynamics of such a one-dimensional WPE in a sinusoidal potential. We find steady states of the system that correspond to a stationary WPE as well as a rich array of unsteady motions such as back-and-forth oscillating walkers, runaway oscillating walkers and various types of irregular walkers. In the parameter space formed by the dimensionless parameters of the applied sinusoidal potential, we observe patterns of alternating unsteady behaviors suggesting interference effects. Additionally, in certain regions of the parameter space, we also identify multistability in the particle's long-term behavior that depends on the initial conditions. We make analogies between the identified behaviors in the WPE system and Bragg's reflection of light as well as electron motion in crystals. 

\end{abstract}

\maketitle

\begin{quotation}
Droplets of oil can bounce and walk on a vertically vibrating bath of the same oil via a coupling between the droplet and its self-generated wave field. This gives rise to a moving wave-particle entity (WPE) on the free surface of the liquid. Such WPEs have been shown to exhibit several features that are typically associated with the quantum realm. In this work, we use a simple theoretical model that takes the form of a Lorenz-like system to investigate the dynamics of a one-dimensional WPE in a sinusoidal potential. We find a rich array of dynamical behaviors in the parameter space which we explore in detail. We also make connections between the dynamics of the WPE and the transport of light and quantum particles in crystals.   
\end{quotation}

\section{Introduction}

A millimeter-sized classical wave-particle entity (WPE) can be realized in the form of a droplet walking on the free surface of a vertically vibrating liquid bath~\cite{Couder2005,Couder2005WalkingDroplets,superwalker}. The walking droplet, also known as a walker or superwalker~\citep{superwalker,superwalkernumerical}, undergoes periodic vertical bouncing motion resonant with the bath driving. Upon each bounce, the droplet generates a slowly decaying localized standing wave on the free surface of the liquid. The interaction of the droplet with its underlying wave field gives rise to horizontal walking motion. For large amplitudes of bath vibrations, the droplet-generated waves decay very slowly in time and the droplet's motion is not only influenced by the most recently generated wave but also by the waves generated in distant past, giving rise to path memory in the system and making the dynamics non-Markovian. In the high-memory regime, walking droplets have been shown to exhibit several hydrodynamic quantum analogs. Examples include quantization of orbits in rotating frames~\citep{Fort17515,harris_bush_2014,Oza2014} and confining potentials~\citep{Perrard2014b,Perrard2014a,labousse2016,PhysRevE.103.053110}, Zeeman splitting in rotating frames~\citep{Zeeman}, surreal Bohmian trajectories~\citep{Surreal2022}, analogs of spin states~\citep{spinstates2018,spinstate2021} and spin systems~\citep{Saenz2021}, tunneling across submerged barriers~\citep{Eddi2009,tunnelingnachbin,tunneling2020}, wave-like statistics in confined geometries~\citep{PhysRevE.88.011001,Giletconfined2016,Saenz2017,Cristea,durey_milewski_wang_2020}, a hydrodynamic analog of Friedel oscillations~\citep{Friedal} and analogues of quantum optics such as hydrodynamic superradiance~\citep{Superradiance2021,Frumkin2021} and Hong-Ou-Mandel-like correlations~\citep{ValaniHOM}. Recent work has also focused on the development of walker-inspired quantum theories~\citep{Dagan2020hqft,Durey2020hqft,mechanalog2020}. A detailed review of hydrodynamic quantum analogs for walking droplets is provided by \citet{Bush2015} and \citet{Bush_2020}.

Simulations of walking-droplet inspired WPEs in a generalized pilot-wave framework have shown a plethora of rich behaviors in the horizontal walking dynamics. The generalized pilot-wave framework~\citep{Bush2015} is a theoretical framework that extends models of walking droplets beyond experimentally achievable regimes. It has allowed for exploration of a broader class of dynamical systems and the discovery of new quantum analogs~\citep{Bush_2020}. In free space with two horizontal dimensions, the typical motion is a steady state where the WPE moves in a straight line at constant speed. However, it has been shown that the steady walking state can become unstable in the high-memory regime and more complex unsteady motions arise such as stable circular orbits, wobbling orbits, precessing orbits and erratic run-and-tumble-like motion~\citep{Bacot2019,Hubert2019,Durey2D2021}. In free space with only one horizontal dimension, the WPE in the unsteady regime has been observed in simulations to exhibits back-and-forth oscillations with and without a net drift as well as chaotic walking motion~\citep{Durey2020,ValaniUnsteady}. Moreover, direct connections have been established between the equation of motion of the WPE in $1$D and Lorenz-like dynamical systems~\citep{Durey2020lorenz,ValaniUnsteady,Valani2022lorenz}.
    
Rich behaviors in the dynamics of a WPE also emerge in simulations within a generalized pilot-wave framework with an added external potential. For a $1$D WPE in a linear tilted potential, it has been shown to exhibits anomalous transport behaviors such as absolute negative mobility, differential negative mobility and lock-in regions corresponding to mobility being independent of the applied force~\citep{Valani2022ANM}. In a harmonic potential, both in $1$D and $2$D, the dynamical constraints imposed by the underlying wave field results in periodic, quasiperiodic, as well as chaotic orbits which may be characterized by intermittent transitions between unstable periodic orbits~\citep{phdthesismolacek,Kurianskiharmonic,durey2018,labousse2016b,Tambasco2016,Perrard2014,Hubert2022}. Similar quantization of orbits as well as chaotic trajectories is also realized in other central potentials and circular corrals~\citep{Tambascoorbit,Cristea,Tambasco2016,Montes2021,Budanur2019,durey_milewski_wang_2020}. The interaction of a WPE with a submerged pillar has been shown to give rise to trajectories in the form of a logarithmic spiral~\citep{Harris2018}. Moreover, macroscopic tunneling-like effects have also been observed where the WPE can unpredictably cross a potential barrier~\citep{hubert2017} or intermittently jumps between two confining cavities in an unpredictable way~\citep{nachbin2017}.

Interaction of light and matter with periodic structures can give rise to emergent phenomena. For example, an electron inside a crystalline solid gives rise to a band structure that forms the basis of emergent material properties such as electrical conductivity~\citep{Kittel2004-xl}. Another example is Bragg's diffraction that takes place when light waves are incident on a crystal whose spacing between the atomic layers is comparable to the wavelength of the incident light~\citep{Bragg1913}. \citet{Braggdroplet2019} explored an analog of Bragg's reflection in the system of walking droplets. They experimentally showed that when a walker is confined to move on an annular track with periodically spaced subsurface barriers, the average speed of the walker nearly vanishes when the spacing between the barriers is close to half the wavelength of the walker-generated damped standing waves. Inspired by these observations in the walking-droplet system, in this work, we use an idealized theoretical model to explore the dynamics of a $1$D walking-droplet inspired WPE in a sinusoidal potential.

The paper is organized as follows. In Sec.~\ref{Sec: theory} we provide details of the theoretical model that governs the dynamics of the WPE. In Sec.~\ref{sec: lin stab full model} we find equilibrium states of the system and determine their stability via a linear stability analysis. In Sec.~\ref{Sec: PS space} we explore the rich dynamics in the parameter space of the system, followed by an exploration of multistability in Sec.~\ref{sec: IC}. In Sec.~\ref{Sec: analogy} we explore connections between the behaviors observed in our system and features associated with propagation of light and quantum particles in crystals. We conclude in Sec.~\ref{Sec: conclusion}.

\section{Theoretical Model}\label{Sec: theory}

\begin{figure}
\centering
\includegraphics[width=\columnwidth]{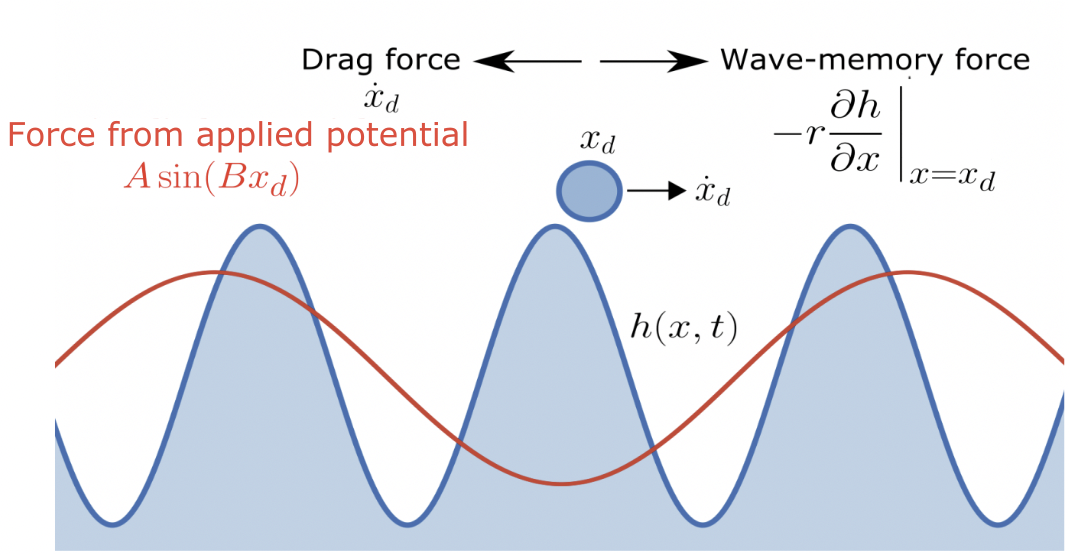}
\caption{Schematic of the system (in dimensionless units) showing a one-dimensional self-propelled WPE (blue) in a sinusoidal potential (red). A particle of dimensionless mass $1/\sigma$ located at $x_d$ and moving horizontally with velocity $\dot{x}_d$ experiences a propulsion force,$-r \partial h/\partial x {|}_{x=x_d}$, from its self-generated wave field $h(x,t)$ (blue filled area), and an effective drag force, $-\dot{x}_d$. The underlying wave field $h(x,t)$ is a superposition of the individual waves of spatial form $W(x)$ and decaying exponentially in time, that are generated by the particle continuously along its trajectory. An external sinusoidal potential results in the particle experiencing an additional force $A\,\sin(B x_d)$.}
\label{Fig: schematic}
\end{figure}

\begin{figure*}
\centering
\includegraphics[width=2\columnwidth]{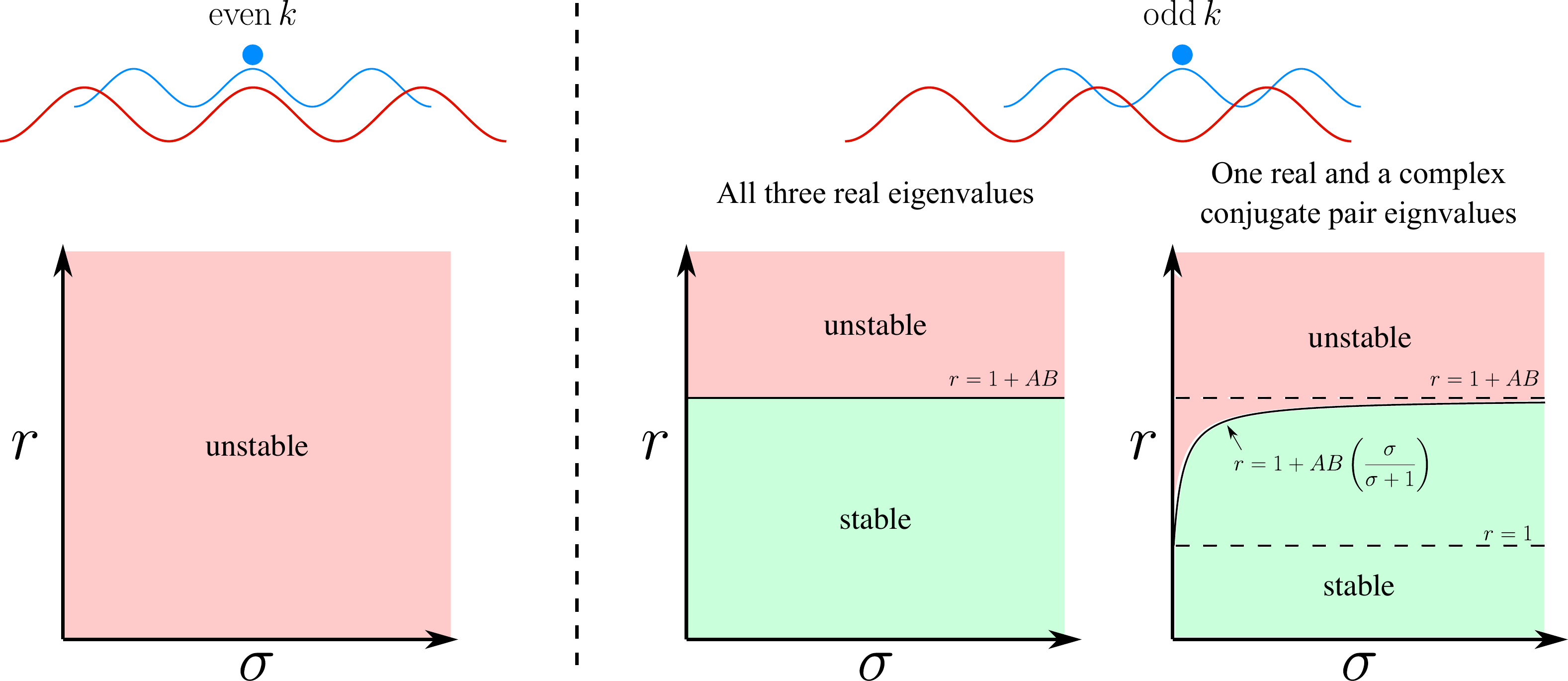}
\caption{Summary of linear stability analysis of the stationary states of a WPE in a sinusoidal potential. For even $k$, corresponding to the particle equilibria at the peaks of the applied sinusoidal potential, the stationary state is always unstable to small perturbations. For odd $k$, corresponding to the particle equilibria at the troughs of the applied sinusoidal potential, the boundary separating stable and unstable regions depends on the nature of the eigenvalues.}
\label{Fig:stabilitysummary}
\end{figure*}

As shown schematically in Fig.~\ref{Fig: schematic}, consider a particle (droplet) located at position $x_d$ and moving horizontally with velocity $\dot{x}_d$ while continuously generating standing waves with prescribed spatial structure $W(x)$ that decay exponentially in time. \citet{Oza2013} developed a theoretical stroboscopic model that averages over the vertical periodic bouncing motion of the particle and provides a trajectory equation for the horizontal walking dynamics in $2$D. A reduction of this model to describe the dynamics of a $1$D WPE with an added external sinusoidal potential is given by the following integro-differential equation of motion:~\citep{Oza2013,ValaniUnsteady} 
\begin{align}\label{Eq: dimensional eq}
    m\ddot{x}_d+D\dot{x}_d =& 
    \frac{F}{T_F}\int_{-\infty}^t f(k_F(x_d(t)-x_d(s)))\,\text{e}^{-(t-s)/(T_F \text{Me})}\,\text{d}s\\ \nonumber
    &+\bar{A}\sin(\bar{k} x_d). 
\end{align}
The left-hand-side of Eq.~\eqref{Eq: dimensional eq} is composed of an inertial term $m \ddot{x}_d$ and an effective drag term $D \dot{x}_d$, where the overdot denotes a time derivative. The first term on the right-hand-side of the equation captures the forcing on the droplet by the underlying wave field. This force is proportional to the gradient of the underlying wave field. The wave field $h(x,t)$ is calculated through integration of the individual wave forms $W(x)$ that are continuously generated by the particle along its trajectory and decay exponentially in time. The function $f(x)=-W'(x)$ is the negative  gradient of the wave form created by the particle. The second term on the right-hand-side is a force arising from an external sinusoidal potential $U(x)=(\bar{A}/\bar{k})\cos(\bar{k}x)$. The parameters in the equation are as follows: $m$ is the droplet mass, $D$ is the effective time-averaged drag coefficient comprising of aerodynamic drag on the droplet during flight and momentum loss during contact with the bath, $k_F=2\pi/\lambda_F$ is the Faraday wavenumber with $\lambda_F$ the Faraday wavelength, $F=m g A_m k_F$ is a nonnegative wave-memory force coefficient with $g$ as the gravitational acceleration and $A_m$ as the amplitude of surface waves, $\text{Me}$ is the memory parameter that describes the proximity to the Faraday instability and $T_F$ is the Faraday period i.e. the period of droplet-generated standing waves and also the bouncing period of the walking droplet. We refer the interested reader to \citet{Oza2013} for more details and explicit expressions for these parameters. The force from applied sinusoidal potential has parameters: $\bar{A}$ denoting the sinusoidal force coefficient, and $\bar{k}=2\pi/\bar{\lambda}$ denoting the wavenumber associated with the wavelength $\bar{\lambda}$ of the applied potential. 
Non-dimensionalizing the equation using $t'={t}/{T_F \text{Me}}$ and $x'=k_F x$, and dropping the primes on the dimensionless variables we get~\citep{Oza2013},
\begin{align}\label{Eq: dimless eq}
\frac{1}{\sigma}\ddot{x}_d+\dot{x}_d=r\int_{-\infty}^t f(x_d(t)-x_d(s))\,\text{e}^{-(t-s)}\,\text{d}s
    +A\sin(Bx_d).
\end{align}
Here, ${\sigma}^{-1}={m}/{DT_F\text{Me}}$ and $r={F \text{Me}^2 T_F k_F}/{D}$ denote the non-dimensional mass and wave-memory force coefficient, respectively. Also, $A={\bar{A}\text{Me}T_Fk_F}/{D}$ is a nonnegative dimensionless force coefficient corresponding to the external sinusoidal potential, while $B={\bar{k}}/{k_F}=\lambda_F/\bar{\lambda}$, denotes the ratio of the wavelength of the droplet-generated waves to the wavelength of the external sinusoidal potential.

We choose a simple sinusoidal particle-generated wave form such that $W(x)=\cos(x)$ and $f(x)=\sin(x)$, which allows us to transform the integro-differential equation of motion in \eqref{Eq: dimless eq} into a Lorenz-like system of ordinary differential equations (ODEs) giving us~\cite{phdthesismolacek,Durey2020lorenz,ValaniUnsteady,Valani2022lorenz} (see \citet{Valani2022lorenz} for a derivation)
\begin{equation}
    \label{lorenz}
    \begin{split}
    \dot{x}_d&=X \\
    \dot{X}&=\sigma(Y-X+A\sin(Bx_d)), \\
    \dot{Y}&=-XZ+r X-Y, \\
    \dot{Z}&=XY-Z.
  \end{split}
\end{equation}

 These ODEs are the Lorenz equations~\citep{Lorenz1963} with an added equation, $\dot{x}_d=X$, and an added feedback term, $\sigma A\sin(Bx_d)$, in the $\dot{X}$ equation in \eqref{lorenz}. Here, $X=\dot{x}_d$ is the droplet's velocity, $Y=r \int_{-\infty}^{t} \sin(x_d(t)-x_d(s))\,\text{e}^{-(t-s)}\,\text{d}s$ is the wave-memory force and $Z=r -r\int_{-\infty}^{t} \cos(x_d(t)-x_d(s))\,\text{e}^{-(t-s)}\,\text{d}s$, is also related to the wave-memory forcing. For the simulations presented in this paper, the system of ODEs in Eq.~\eqref{lorenz} is solved in MATLAB using the inbuilt solver ode45. We note that the initial conditions for these simulations correspond to fixed points of the dynamical system; however, there is an intrinsic perturbation provided by numerical round-off errors in MATLAB resulting in unsteady dynamics when the fixed points are unstable. Moreover, for all results presented in this paper, the system was simulated till $t=2000$.


\section{Equilibrium Solutions and linear stability analysis}\label{sec: lin stab full model}

We start by finding equilibrium solutions of the dynamical system. Setting the time derivatives to zero in Eq.~\eqref{lorenz} results in the following infinite set of equilibrium solutions: 
\begin{equation*}
    x_d=\frac{k\pi}{B},\:X=0,\:Y=0\:\:\text{and}\:\:Z=0,
\end{equation*}
where $k\in\mathbb{Z}$. This corresponds to a stationary WPE located at the peaks and troughs of the applied sinusoidal potential. We note that for this WPE in free space (i.e. in the absence of the external sinusoidal potential), the dynamical system also has an equilibrium solution corresponding to a steadily moving WPE~\citep{Valani2022lorenz,ValaniUnsteady}. This solution is 
\begin{equation*}
x_d=x_0+Xt,\:X=\pm\sqrt{r-1},\:Y=\pm\sqrt{r-1}\:\:\text{and}\:\:Z=r-1,
\end{equation*}
where $x_0$ is any real number and it corresponds to translational invariance in the absence of the external potential. 

To determine the stability of the stationary solution in the presence of the sinusoidal potential, we perform a linear stability analysis by applying a small perturbation to the equilibrium solutions as follows~\citep{strogatz}: $(x_d,X,Y,Z)=(k\pi/B,0,0,0)+\epsilon(x_{d1},X_1,Y_1,Z_1)$, where $\epsilon>0$ is a small perturbation parameter. Substituting this in Eq.~\eqref{lorenz} and comparing terms of $O(\epsilon)$, we obtain the following linear system that governs the evolution of perturbations:
\begin{gather*}
 \begin{bmatrix} 
 \dot{x}_{d1} \\
 \dot{X}_1 \\
 \dot{Y}_1 \\
 \dot{Z}_1 
 \end{bmatrix}
 =
  \begin{bmatrix}
0 & 1 & 0 & 0 \\
(-1)^k AB\sigma & -\sigma & \sigma & 0 \\
0 & r & -1 & 0 \\
0 & 0 & 0 & -1
 \end{bmatrix}
  \begin{bmatrix}
  x_{d1}\\
  {X}_1 \\
 {Y}_1 \\
 {Z}_1 
 \end{bmatrix}.
\end{gather*}
The linear stability is determined by the eigenvalues of the right-hand-side matrix. This results in the following characteristic polynomial equation to be solved for the eigenvalues $\lambda$ which determines the growth rate of perturbations:
\begin{align}\label{eq: polynomial}
    \left( \lambda+1 \right) \Big[ \lambda^3+(\sigma+1)\lambda^2&+\sigma\left(1-r-(-1)^k AB\right)\lambda\\ \nonumber
    & -(-1)^k AB\sigma \Big]=0
\end{align}
To obtain non-trivial eigenvalues of the characteristic polynomial in Eq.~\eqref{eq: polynomial} requires solving a cubic equation. However, by invoking Descartes' rule of sign, one can deduce information about the nature of the real eigenvalues of this cubic equation without explicitly solving the cubic equation. 

For even $k$, we use Descartes' rule of sign to deduce the existence of one positive real root of the cubic equation. This existence of a real positive eigenvalue informs us that for even $k$, which corresponds to the particle equilibria at the peaks of the applied potential, the stationary state is always unstable to small perturbations.

%
%
%
%
%
%

For odd $k$, corresponding to the particle equilibria at the troughs of the applied sinusoidal potential, it follows from Descartes' sign rule that: (i) if $r<1+AB$, then there are zero real positive eigenvalues and either three or one negative real eigenvalues. So either all the eigenvalues are real and negative or one of them is a negative real number and the two others form a complex conjugate pair. (ii) If $r>1+AB$, we find that there are either two or zero real positive eigenvalues and one real negative eigenvalue. So either we have two real positive and one real negative eigenvalues or one real negative eigenvalue and a pair of complex conjugate eigenvalues. Thus, if all the eigenvalues are real, the stationary state is stable to small perturbations for $r<1+AB$ and unstable to small perturbations for $r>1+AB$. Alternatively, when a complex conjugate pair of eigenvalues exists, one eigenvalue is always real and negative for all $r$ and hence the stability of the stationary solution will depend on the nature of the complex conjugate pair. The criteria for existence of such a complex conjugate pair can be deduced from the discriminant $\Delta$ of the cubic in Eq.~\eqref{eq: polynomial} given by
\begin{align*}
    \Delta&=\alpha^2_1 \alpha^2_2 - 4 \alpha^3_2 - 4\alpha^3_1 \alpha_3 - 27 \alpha^2_3 + 18 \alpha_1 \alpha_2 \alpha_3,
\end{align*}
with
\begin{align*}
    \alpha_1&=\sigma + 1, \\
    \alpha_2&=\sigma (1-r+AB), \\
    \alpha_3&=AB\sigma.
\end{align*}
There will be one real and a pair of complex conjugate roots when $\Delta<0$ and all real roots when $\Delta>0$.

To understand the nature of the complex conjugate eigenvalues for odd $k$, we determine the boundary that corresponds to complex conjugate eigenvalues changing the sign of its real part i.e. the boundary corresponding to $\text{Re}(\lambda)=0$. Substituting $\lambda=i\omega$ into the cubic polynomial in Eq.~\eqref{eq: polynomial} and setting the real and imaginary parts to zero gives us the following set of solutions:

%
%
%
\begin{equation}\label{stability boundary}
\begin{split}
    r=r_c=1+AB\left(\frac{\sigma}{\sigma+1}\right),
    \end{split}
\end{equation}
and
\begin{equation}\label{small oscillations}
    \omega^2=AB \left(\frac{\sigma}{\sigma+1}\right).
\end{equation}


The curve given by Eq.~\eqref{stability boundary} gives the instability boundary in the $(\sigma, r)$ parameter space. Below this boundary the complex conjugate eigenvalues have $\text{Re}(\lambda)<0$ and above this boundary we have $\text{Re}(\lambda)>0$. Note that the boundary only depends on the product $AB$ of the parameters of the external force. Thus, we find that when complex conjugate eigenvalues exist, for $r<r_c$ the real part of the complex conjugate pair is negative which corresponds to the stationary state being stable to small perturbations, and for $r>r_c$ the positive real part of the complex conjugate eigenvalues results in growth in oscillations and the equilibrium solution is unstable to small perturbations. Moreover, at the onset of the instability just above the stability boundary, the frequency of small oscillations is given by Eq.~\eqref{small oscillations}. The stability analysis is visually summarized in Fig.~\ref{Fig:stabilitysummary}.

\begin{figure*}
\centering
\includegraphics[width=2\columnwidth]{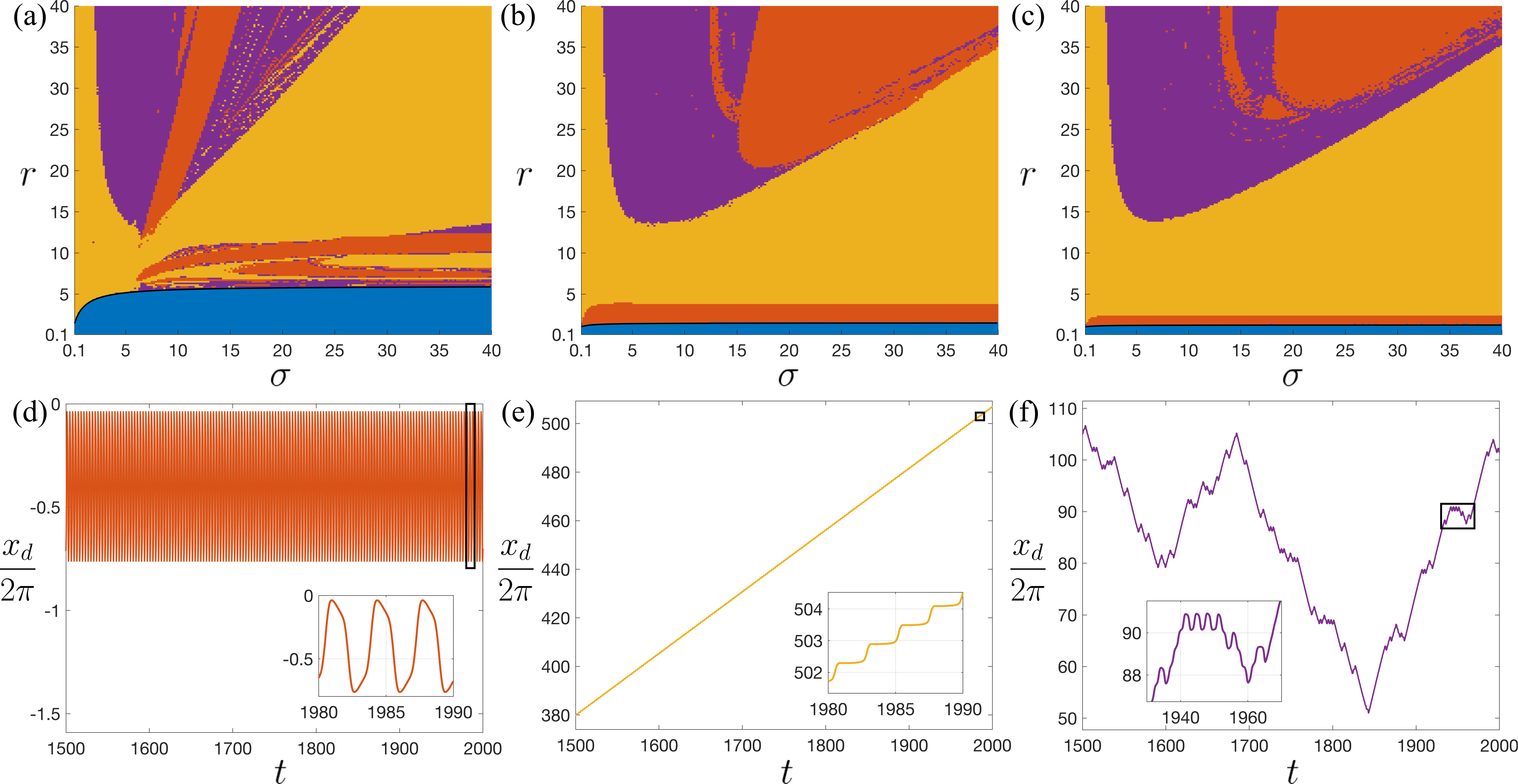}
\caption{Dynamical behaviors of a WPE in a sinusoidal potential. Different observed behaviors in the $(\sigma, r)$ parameter space for fixed (a) $(A,B)=(1,5)$, (b) $(A,B)=(1,0.5)$ and (c) $(A,B)=(0.5,0.5)$. In these plots, blue corresponds to stationary states, red indicates back-and-forth oscillating walkers, yellow indicates runaway oscillating walkers and purple indicates irregular chaotic walkers. Example space-time trajectories of the three qualitatively different unsteady behaviors are shown as (d) back-and-forth oscillating walkers $(A,B,\sigma,r)=(1,5,8,15)$, (e) runaway oscillating walkers $(A,B,\sigma,r)=(1,5,20,15)$ and (f) irregular walkers $(A,B,\sigma,r)=(1,5,4.5,19.5)$. The black curve in panels (a)-(c) shows the linear stability curve of the stationary state calculated from Eq.~\eqref{stability boundary}. The initial conditions were fixed to $(x_d(0),X(0),Y(0),Z(0))=(\pi/B,0,0,0)$ which corresponds to a stationary state.}
\label{Fig:dynamicsPS}
\end{figure*}

\begin{figure*}
\centering
\includegraphics[width=1.3\columnwidth]{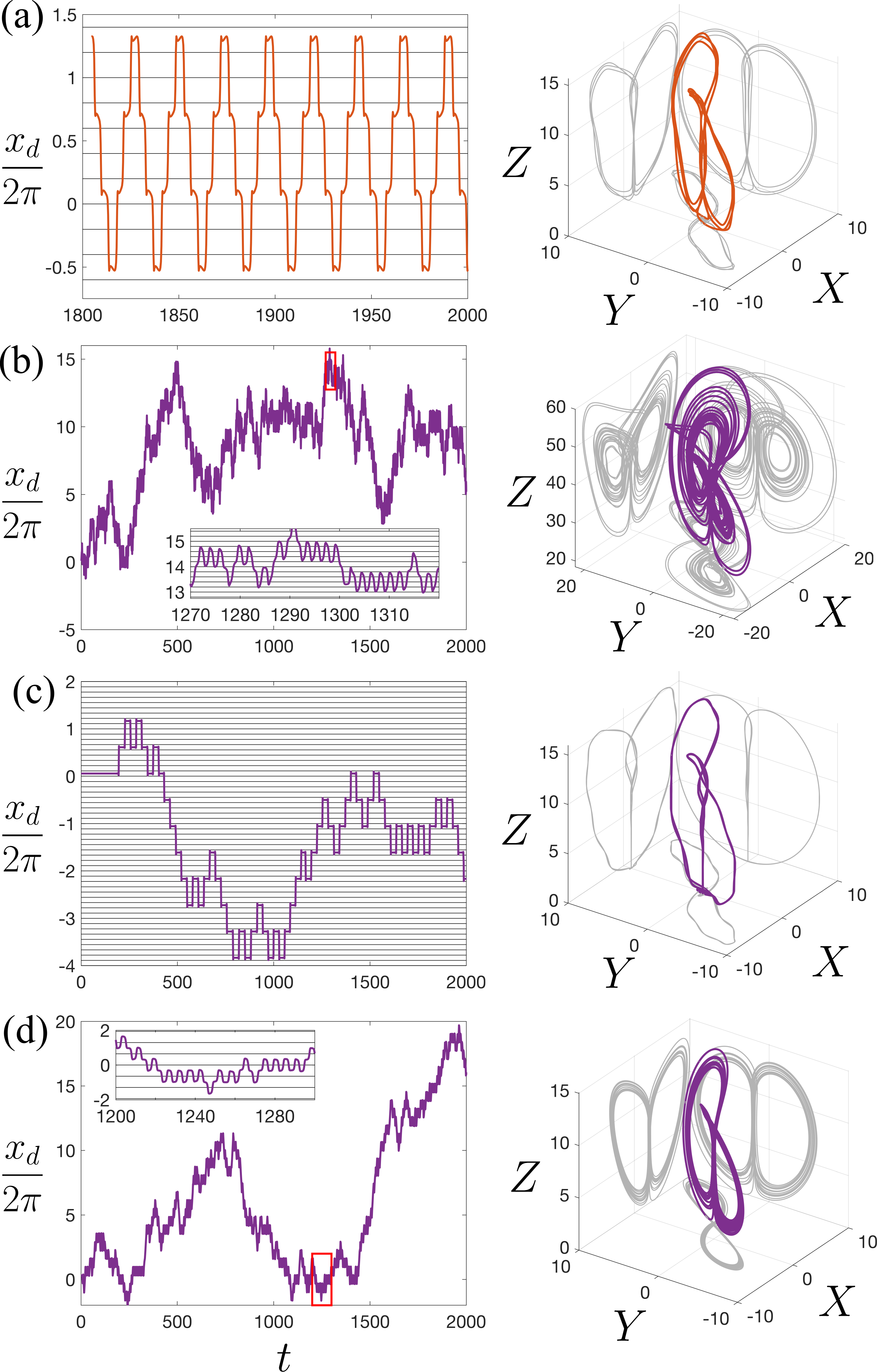}
\caption{Intricate unsteady dynamical behaviors of a WPE in a sinusoidal potential. Space-time trajectories and phase-space dynamics are shown for (a) a back-and-forth oscillating walker with complex oscillations, $(A,B,\sigma,r)=(1,5,10,10)$, (b) an irregular walker with a Lorenz-like chaotic attractor, $(A,B,\sigma,r)=(1,5,10,40)$, (c) an irregular walker that remains stationary in a local potential well for a fixed duration and then jumps unpredictably to the left or right, $(A,B,\sigma,r)=(1.05,9,10,10)$, and (d) an irregular walker that performs an unpredictable number of back-and-forth oscillations between consecutive wells of the applied potential before jumping unpredictably to neighboring pair of wells, $(A,B,\sigma,r)=(0.6,1.55,10,10)$. The black horizontal lines represent the maxima of the applied sinusoidal potential. The initial conditions for these simulations were $(x_d(0),X(0),Y(0),Z(0))=(\pi/B,0,0,0)$.}
\label{Fig: intricatebehav}
\end{figure*}

\section{Parameter space exploration}\label{Sec: PS space}

\begin{figure*}
\centering
\includegraphics[width=2\columnwidth]{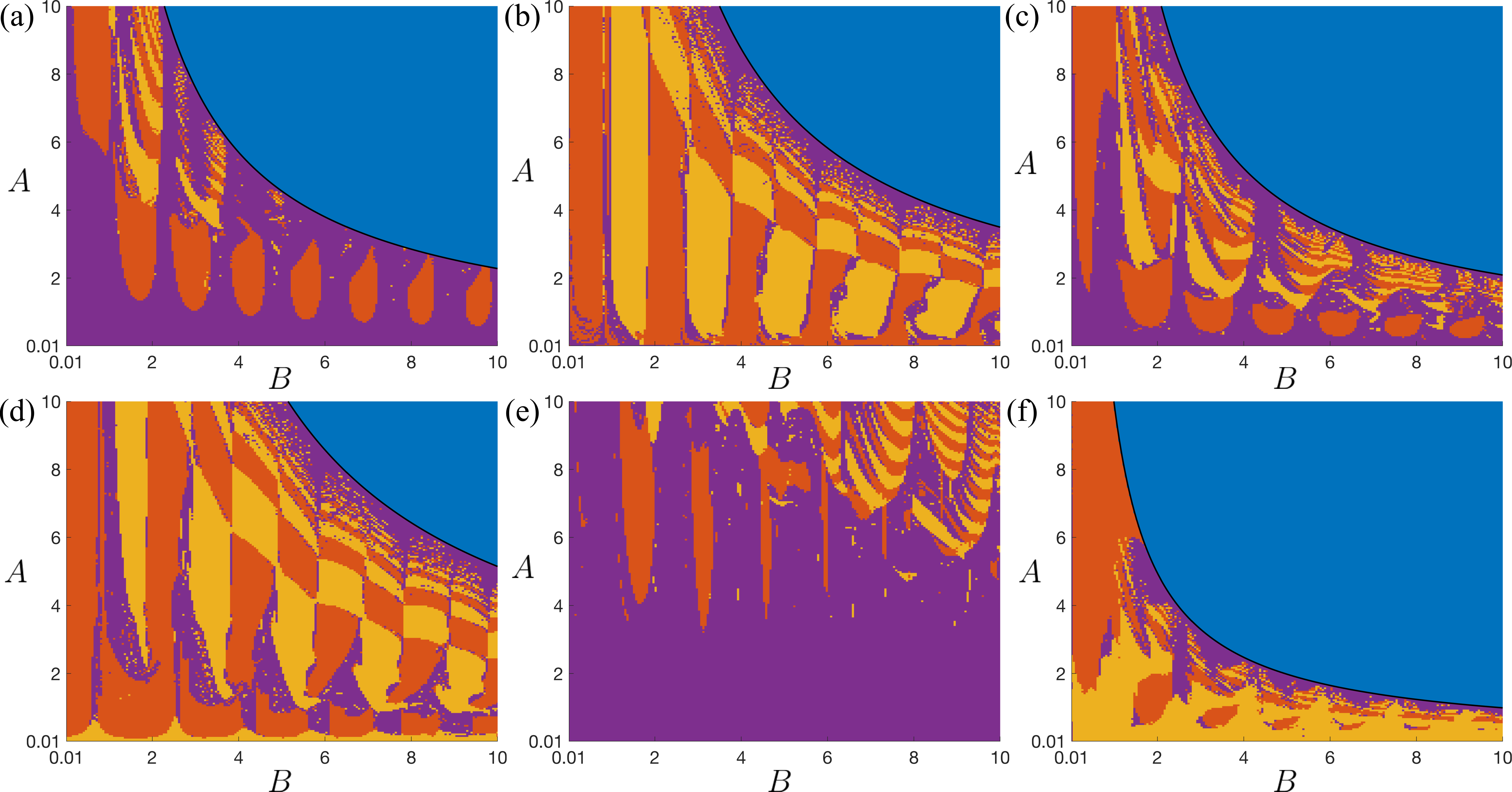}
\caption{Interference effects in the dynamical behaviors of a WPE in a sinusoidal potential. Different dynamical behaviors in the $(A,B)$ parameter space for fixed (a) $(\sigma,r)=(5,20)$, (b) $(\sigma,r)=(35,35)$, (c) $(\sigma,r)=(10,20)$, (d) $(\sigma,r)=(20,50)$, (e) $(\sigma,r)=(10,100)$ and (f) $(\sigma,r)=(10,10)$. Interference effects are observed as alternating unsteady behavior patterns in these parameter space plots. Here, blue corresponds to stationary states, red indicates back-and-forth oscillating walkers, yellow indicates runaway oscillating walkers and purple indicates irregular walkers. The black curve shows the linear stability curve of the stationary state calculated from Eq.~\eqref{stability boundary}. The initial conditions were taken as $(x_d(0), X(0), Y(0), Z(0))= (\pi/{B},0,0,0)$.}
\label{Fig: AB space}
\end{figure*}

To explore the range of dynamical behaviors exhibited by the WPE, we simulate space-time particle trajectories for different values of the system parameters $\sigma$, $r$, $A$ and $B$.

We start by exploring the dynamical behaviors of the WPE in the parameter space formed by $\sigma$ and $r$. Different types of observed behaviors in this parameter space are shown in Fig.~\ref{Fig:dynamicsPS}(a)-(c) for three different sets of $(A,B)$ values. We classify the observed long-time behaviors at the end of our simulations into four types: (i) stationary states (blue) where the particle remains stationary with no horizontal motion, (ii) back-and-forth oscillating walkers (red) where the particle gets trapped and undergoes back-and-forth motion about a fixed location, (iii) runaway oscillating walkers (yellow) where the particle undergoes oscillatory motion along with a net drift, and (iv) irregular walkers (purple) where the particle exhibits a chaotic diffusive motion. Example trajectories corresponding to unsteady behaviors (ii)-(iv) are shown in Fig.~\ref{Fig:dynamicsPS}(d)-(f), respectively. We note that similar types of unsteady behaviors are also observed for the WPE without the external sinusoidal potential in the parameter regime of large $\sigma$ and $r$, where the strong interaction of the particle with its self-generated sinusoidal wave field results in these unsteady states~\citep{ValaniUnsteady}. Moreover, similar trajectories have also been observed in a generalized Lorenz system where the motion of a fictitious particle was considered in a periodic lattice of truncated parabolae~\citep{Festa_2002}. 

Regions occupied by each of the behaviors in the $(\sigma,r)$ parameter space vary with $A$ and $B$ values. The region occupied by stationary states (blue) can be deduced from linear stability analysis and it occupies the region below the linear stability curve (black curve) in Fig.~\ref{Fig:dynamicsPS}(a)-(c) as calculated from Eq.~\eqref{stability boundary}. Above this curve, the motion of the particle becomes unsteady. For $A=1$ and $B=5$, as shown in Fig.~\ref{Fig:dynamicsPS}(a), we find that a large region of the parameter space above the linear stability curve is occupied by runaway oscillatory walkers (yellow). We also have a lobe-shaped region consisting mainly of irregular walkers (purple) with a band of back-and-forth oscillating walkers (red) and a scattered mixture of all three kinds of unsteady behaviors. 
Just above the linear stability curve for $\sigma \gtrsim 6$, we also obtain a mixture of all three unsteady behaviors. Behaviors realized in the parameter space for $(A,B)=(1,0.5)$ and $(0.5,0.5)$ as shown in Fig.~\ref{Fig:dynamicsPS}(b) and (c), respectively, reveal similar lobe-shaped regions consisting of irregular walkers and back-and-forth oscillating walkers. Moreover, for these latter parameter space plots, we observe only back-and-forth oscillating walkers just above the linear stability curve for stationary states. We note that the behaviors observed in some regions of the parameter space are sensitive to initial conditions and we will further explore this aspect in Sec.~\ref{sec: IC}. 


We have also observed more intricate types of unsteady behaviors in the parameter space. A few examples of such behaviors showing space-time trajectories along with the corresponding $(X,Y,Z)$ phase-space trajectories are depicted in Fig.~\ref{Fig: intricatebehav}. For example, both runaway oscillating walkers and back-and-forth oscillating walkers can have complex oscillations in one-period rather than simple sinusoidal oscillations. An example of a space-time trajectory for such oscillations and the corresponding phase-space trajectory are shown in Fig.~\ref{Fig: intricatebehav}(a) for a back-and-forth oscillating walker. Here, we note that the amplitude of oscillations of the back-and-forth walker is on the order of the wavelength of its self-generated waves. Moreover, the turning points of oscillations are correlated with the wells of the external sinusoidal potential (black horizontal lines represent maxima of the external sinusoidal potential in Fig.~\ref{Fig: intricatebehav}). Chaotic irregular walkers also have several distinctive features. A typical irregular walker is shown in Fig.~\ref{Fig: intricatebehav}(b). This irregular walker shows no apparent order in its trajectory, but its phase-space dynamics reveals a Lorenz-like strange attractor. Moreover, as shown in the inset of Fig.~\ref{Fig: intricatebehav}(b), the walker performs back-and-forth oscillatory motion at small time scales. The amplitude of oscillations are of the order of the wavelength of its self-generated waves and for this walker the turning points are typically associated with the maxima of the applied potential. We have also observed other types of irregular walkers where the WPE either remains stationary in a well of the applied potential and then intermittently jumps between neighboring wells (see Fig.~\ref{Fig: intricatebehav}(c)) or performs back-and-forth periodic oscillations between a pair of consecutive wells of the applied potential before hopping unpredictably to another consecutive pair of wells (Fig.~\ref{Fig: intricatebehav}(d)). For the former, as shown in Fig.~\ref{Fig: intricatebehav}(c), the time spent in each well is almost constant and the only unpredictability is in selecting to jump to the left or right neighbouring well. Moreover, the trajectory also shows hints of pseudolaminar chaotic diffusion~\citep{Muller2023}. For the latter, as shown in Fig.~\ref{Fig: intricatebehav}(d), the particle performs an unpredictable number of oscillations between consecutive potential wells before jumping to a neighboring pair of wells also in an unpredictable way.

To understand the effects of the external potential on the particle dynamics, we also explore how the dynamical behavior changes in the parameter space formed by $A$ and $B$. Intricate structures of alternating unsteady behaviors are shown for six different sets of $(\sigma, r)$ values in Fig.~\ref{Fig: AB space}. For example in Fig.~\ref{Fig: AB space}(b) where $(\sigma,r)=(35,35)$, we observe a checkerboard pattern of alternating behaviors between back-and-forth oscillatory walkers and runaway oscillatory walkers. In such a parameter-space, a horizontal transect at constant force coefficient (fixed $A$) or a vertical transect at constant wavelength ratio (fixed $B$) would lead to alternating dynamical behaviors suggesting interference effects between the dynamics of the  WPE and the external sinusoidal potential. We discuss the implications of such interference effects in the context of an analogy with Bragg's reflection in Sec.~\ref{Sec: Bragg}.


\begin{figure}
\centering
\includegraphics[width=\columnwidth]{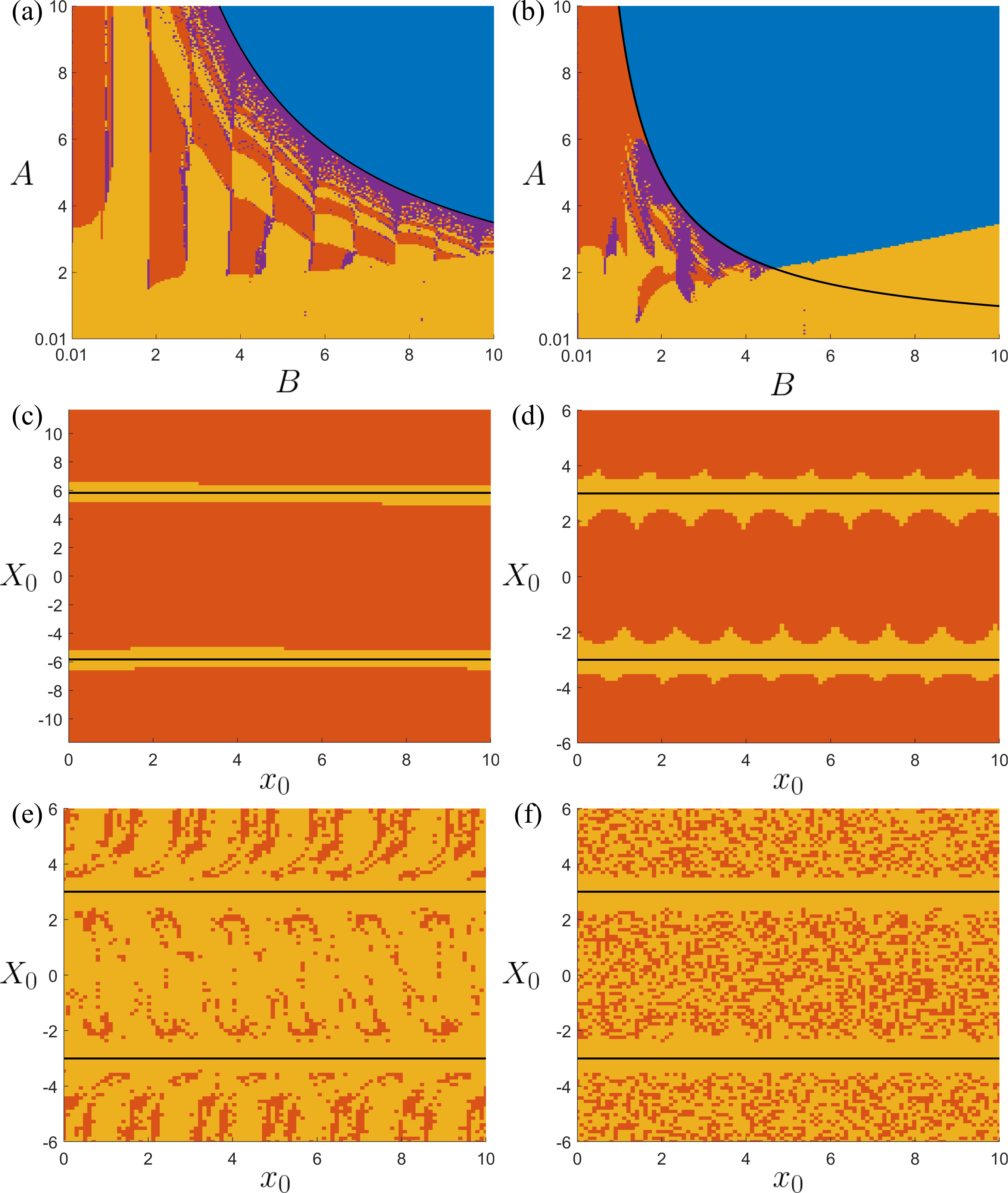}
\caption{Multistability in the dynamics of a WPE in a sinusoidal potential. Dynamical behaviors in the $(A,B)$ parameter space for initial conditions corresponding to a free steady walker: $(x_d(0), X(0), Y(0), Z(0))= (\pi/B,\sqrt{r-1},\sqrt{r-1},r-1)$ and fixed parameters (a) $(\sigma,r)=(35,35)$ and (b) ($\sigma,r)=(10,10)$. Here, blue corresponds to stationary states, red indicates back-and-forth oscillating walkers, yellow indicates runaway oscillating walkers and purple indicates irregular walkers. Basin of attraction of dynamical behaviors in the parameter-space formed by initial position $x_0$ and initial velocity $X_0$ for (c) $(A,B,\sigma,r)=(2,0.5,35,35)$, (d) $(A,B,\sigma,r)=(1,5,10,10)$, (e) $(A,B,\sigma,r)=(1,3.9,10,10)$ and (f) $(A,B,\sigma,r)=(1,5.9,10,10)$. The black horizontal solid lines correspond to the free walking velocities $\pm\sqrt{r-1}$. For panels (c)-(f), the initial conditions were chosen to be $(x_d(0), X(0), Y(0), Z(0))= (x_0,X_0,X_0,X_0^2)$.}
\label{fig: multistability}
\end{figure}

\subsection{Low-memory regime}\label{Sec: low mem}
The theoretical model implemented in this paper arises from the integro-differential equation~\eqref{Eq: dimless eq} which takes into account all the past waves generated by the particle. To understand the role of memory in giving rise to complex dynamical behaviors, we analyze the model in the low-memory regime. We assume that the particle is only influenced by its most recent wave and neglect all other past waves, which is a reasonable assumption at low memory where the waves generated by the particle decay quickly. This results in the following reduced equation of motion~\citep{Protiere2006,Oza2013} (see Appendix~\ref{App: low memory} for a derivation)

\begin{equation}\label{eq: lowmemode}
    \frac{1}{\sigma}\ddot{x}_d+\dot{x}_d=\frac{r}{e}\sin(\dot{x}_d)+A\sin(B x_d),
\end{equation}
which leads to the following system of nonlinear ODEs,
\begin{equation}
\begin{split}
    \label{eq: lowmem}
    \dot{x}_d&=v \\
    \dot{v}&=\sigma\left(\frac{r}{e}\sin(v)+A\sin(B x_d)-v\right).
    \end{split}
\end{equation}
Since the above continuous dynamical system is autonomous and has two-dimensions, it can be concluded from Poincaré–Bendixson theorem that this system cannot exhibit chaotic behavior~\citep{strogatz}. Hence, all the complex irregular walking behaviors in the dynamics of a WPE in a sinusoidal potential arise from path-memory in the system. Nevertheless, we briefly explore the dynamics of this low-memory system in Appendix~\ref{App: low memory}.

\section{Multistability}\label{sec: IC}

We have observed the existence of multiple stable states of the WPE which depend on the initial conditions. For example, if the initial conditions of the WPE correspond to a free steady walking state $(x_d(0), X(0), Y(0), Z(0))= (\pi/B,\sqrt{r-1},\sqrt{r-1},r-1)$ instead of a stationary initial state $(x_d(0), X(0), Y(0), Z(0))= (\pi/B,0,0,0)$, then this leads to a change in the type of behavior observed at long times in some regions of the parameter space. Comparing the $(A,B)$ parameter space behaviors for fixed $(\sigma,r)=(35,35)$ arising in the free-walking initial state, Fig.~\ref{fig: multistability}(a), and the stationary initial state, Fig.~\ref{Fig: AB space}(b), we note that for low values of $A$, the checkerboard pattern present in Fig.~\ref{Fig: AB space}(b) is replaced by runaway oscillating walkers. A similar change is also observed between Fig.~\ref{Fig: AB space}(f) and Fig.~\ref{fig: multistability}(b) for a different choice of $\sigma$ and $r$ values. The presence of runaway oscillating walkers above the linear stability curve of stationary states (black curve) in Fig.~\ref{fig: multistability}(b) shows that the stationary state coexists with the runaway oscillating walkers. These two different behaviors under identical parameter values indicate the sensitivity of the WPE to initial conditions. To further explore the dependence on initial conditions, we have plotted the basin of attraction showing different dynamical behaviors in the parameter-space formed by initial position and initial velocity as shown in Figs.~\ref{fig: multistability}(c)-(f). As it can be seen in Fig.~\ref{fig: multistability}(c), if the initial velocity of the particle is near the free steady walking speed (black horizontal line) then a runaway oscillating walker is realized, while for most other initial conditions, a back-and-forth oscillating walker is realized. For different combinations of parameter values, these basin of attraction plots in the parameter space of initial conditions show more intricate structures as shown in Figs.~\ref{fig: multistability}(d)-(f). The periodic features in Figs.~\ref{fig: multistability}(d) and (e) are associated with periodicity of the applied potential.

\section{Analogy with propagation of light and quantum particles in crystals}\label{Sec: analogy}


\begin{figure*}
\centering
\hspace*{-0.5cm} 
\includegraphics[width=2\columnwidth]{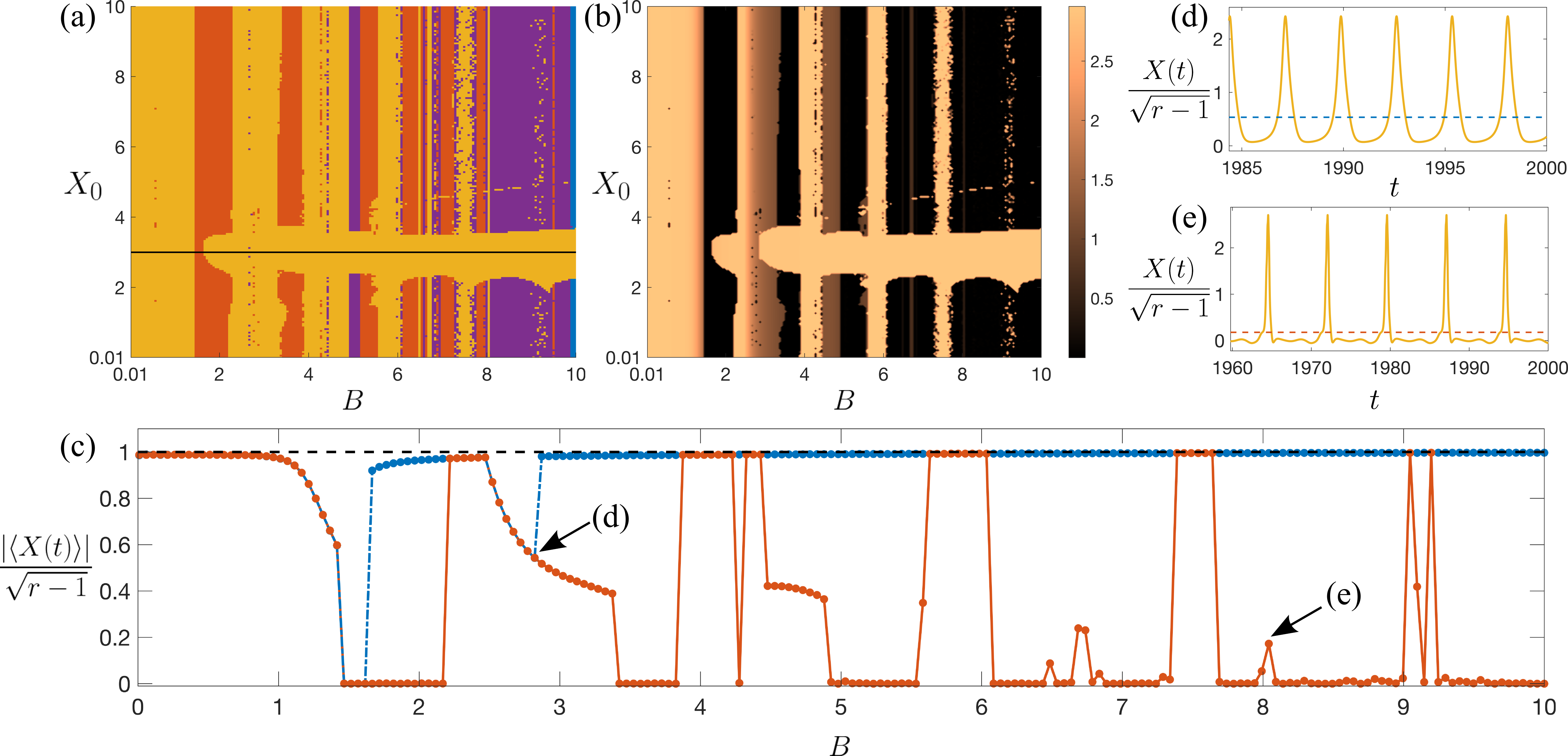}
\caption{Interference effects as a function of the wavelength ratio $B$. (a) Dynamical behaviors in the parameter space formed by the initial velocity $X_0$ and the wavelength ratio $B$. Here, blue corresponds to stationary states, red indicates back-and-forth oscillating walkers, yellow indicates runaway oscillating walkers and purple indicates irregular walkers. The black horizontal line represents the initial velocity corresponding to the free walking speed of the particle, $\sqrt{r-1}$, in the absence of the external sinusoidal potential. (b) Contour plot in the same parameter space with the colorbar showing the absolute average velocity calculated from the velocity time-series, $|\langle X(t) \rangle|=|(1/N)\sum_{i=1}^{N} X(t_i)|$, at the end of simulation time of $2000$. (c) Plot of absolute average velocity normalized by free walking speed $|\langle X(t) \rangle|/\sqrt{r-1}$ against $B$ for $X(0)=\sqrt{r-1}$ (blue dotted-dashed curve) and $X(0)=0$ (red solid curve). Panels (d) and (e) respectively show velocity time series of runaway oscillating walkers at the parameter values corresponding to the first dip in blue curve and a small peak in red curve in panel (c) (see black arrows). The horizontal dashed line represents the average velocity in these panels. Parameter values for these plots are as follows: $A=1, \sigma=10, r=10$ and initial conditions are $(x_d(0),X(0),Y(0), Z(0))=({\pi}/{B},X_0, X_0, X_0^2)$.}
\label{Fig: braggdrop}
\end{figure*}

\subsection{Analogy with Bragg's reflection}\label{Sec: Bragg}

When light waves are incident on a crystal, they transmit through the crystal if there is a large disparity between the wavelength of incident light and the inter-atomic spacing inside the crystal. However, Bragg's diffraction can take place when the wavelength of the incident light is comparable to the inter-atomic spacing inside the crystal, resulting in scattering of the incident light and high intensity reflections. Similar effects are also observed for beams of quantum particles such as electrons or neutrons. An analog of Bragg's reflection has also been reported in the system of walking droplets~\citep{Braggdroplet2019}. By confining a walking droplet to move on an annular track with periodically spaced subsurface barriers, it was found that the average droplet speed nearly vanishes when the periodic subsurface pattern has a characteristic length close to half the Faraday wavelength~\citep{Braggdroplet2019}. Here the authors explained the observed phenomenon by considering the fact that the Faraday instability is altered by the periodic pattern of submerged barriers, which results in generation of small amplitude waves for the walker and hence lower walking speed. This low-speed walker can reflect from subsurface barriers, giving rise to a dip in the plot of the average walking speed as a function of the spacing between the submerged barriers and resulting in a Bragg's reflection-like phenomenon in this system. We can also explore an analog of Bragg's reflection using the simple model presented in this work. However, we note that our model assumes a predefined form of the waves generated by the particle/droplet and hence it does not take into account the changes in the wave amplitude of the droplet-generated waves due to the presence of the underlying periodic potential. Therefore, we do not expect to capture the mechanism of Bragg's reflection-like phenomenon identified by \citet{Braggdroplet2019}. Nevertheless, we have found interference effects even in our simple model as identified in Fig.~\ref{Fig: AB space} and we further discuss them here.

Figures~\ref{Fig: braggdrop}(a) and (b) show the dynamical behaviors and contours of the absolute average velocity of the particle in the parameter space formed by the initial velocity $X_0$ and the wavelength ratio $B$. From Fig.~\ref{Fig: braggdrop}(a) we find that when the initial velocity of the particle is near the free walking speed (black horizontal line), $X_0=\sqrt{r-1}$, runaway oscillating walkers are obtained for all $B$ values shown except a small region of back-and-forth oscillating walkers near $B\approx 1.5$. Moreover, on plotting the absolute average velocity $|\langle X(t) \rangle|$ of the particle as a function of $B$ we find dips in the absolute average velocity (see Fig.~\ref{Fig: braggdrop}(b) and blue dotted-dashed curve in Fig.~\ref{Fig: braggdrop}(c)). The first dip to zero average velocity near $B\approx 1.5$ corresponds to back-and-forth oscillating walkers, while the smaller second dip near $B\approx 3$ is due to the change in nature of oscillations for the runaway oscillating walker resulting in reduced average velocity (see Fig.~\ref{Fig: braggdrop}(d)). If the initial velocity of the particle is not near the free walking speed then we obtain a more complex dependence on $B$ and find several dips in the absolute average velocity. For example, the red solid curve in Fig.~\ref{Fig: braggdrop}(c) corresponding to a stationary initial state, $X_0=0$, shows several dips and complex structure. The small peaks evident in the red solid curve are associated with nonlinear oscillations of runaway oscillating walkers~(for example see Fig.~\ref{Fig: braggdrop}(e)).

The existence of different dynamical behaviors and changes in the absolute average velocity of the particle as a function of $B$ indicates that the interference of the two potentials experienced by the droplet, the droplet-generated wave field potential and the external sinusoidal potential, can affect the transport of the particle across the applied sinusoidal potential. These observation suggest further exploration of these interference effects in connection to analogs with Bragg's reflection, however, they are beyond the scope of the present work.

\subsection{Analogy with electron transport in a crystal}

The motion of electrons in a crystal can give rise to interesting features such as conduction bands and band gaps that form the basis of electrical properties of various materials. \citet{Arpornthip2009} found that for a classical particle in a sinusoidal potential whose energy is allowed to be a complex number, one can get particle motion similar to that of an electron in a crystal. They found two qualitatively different types of particle motion: (i) a random-walk-like hopping motion between neighboring wells of the periodic potential behaving like a quantum particle in an energy gap which undergoes repeated tunneling processes, and (ii) a drifting motion with a constant average speed behaving like a quantum particle in a conduction band undergoing resonant tunneling. The trajectories that we have observed for the WPE are also similar to that observed by \citet{Arpornthip2009} and hence we can make an analogy between the motion of the WPE and that of a quantum particle, e.g. electron, in a crystal. The stationary states and the back-and-forth oscillating walkers are localized states where the particle is either confined in a potential well or oscillating between neighboring potential wells of the periodic potential. These states can be thought of as being analogous to localized states of an electron in a material with large band gap resulting in immobile electrons and the corresponding material being an insulator. Typical runaway oscillating walkers have a constant average speed comparable to its free walking speed and the particle drifts through the periodic potential making these states analogous to delocalized states of electrons in a conduction band. For the irregular walkers, we have seen trajectories (see Fig.~\ref{Fig: intricatebehav}(c,d)) where the walkers spends some time in a potential well or a pair of wells before unpredictably hopping to a neighboring potential wells and hence these states can be thought of as analogous to a quantum particle undergoing a tunneling process between neighboring sites. It would be interesting to explore this analogy in more detail by exploring analogs of solid-state physics in our system of classical WPEs.

\section{Conclusions}\label{Sec: conclusion}

In this paper, we have explored the dynamics of a walking-droplet inspired $1$D classical WPE in a sinusoidal potential. The integro-differential equation of motion conveniently transformed into a system of Lorenz-like ODEs with a sinusoidal feedback term. We found equilibrium stationary states for the particle at peaks and troughs of the applied potential. Using a linear stability analysis we showed that the stationary states at peaks of the applied sinusoidal potential are always unstable while those located at troughs are stable in certain regions of the parameter space~(Fig.~\ref{Fig:stabilitysummary}). Where stationary states are unstable in the parameter space, we observed rich unsteady dynamical behaviors that we explored using numerical simulations~(Figs.~\ref{Fig:dynamicsPS}, \ref{Fig: AB space}). Three qualitatively different types of unsteady behaviors were found: back-and-forth oscillating walkers, runaway oscillating walkers and irregular walkers~(Figs.~\ref{Fig:dynamicsPS}, \ref{Fig: intricatebehav}). We made analogies between these different types of trajectories and the transport of electrons in crystals~(Fig.~\ref{Fig: braggdrop}). 

In addition to rich dynamical behaviors, we also observed interference effects and multistability in the system. Interference effects were realized as alternating unsteady dynamical behaviors as the dimensionless external force coefficient $A$ and the wavelength ratio $B$ were varied~(Fig.~\ref{Fig: AB space}). These interference effects also inspired an analogy with Bragg's reflection where we showed dips in the absolute average particle velocity as a function of the wavelength ratio $B$~(Fig.~\ref{Fig: braggdrop}). We also observed multistability in the system where different initial conditions led to different long-term dynamical behaviors~(Fig.~\ref{fig: multistability}). In the analogy with Bragg's reflection, this resulted in different structures of the dips for different initial velocities of the particle.

The work presented here can be extended in several ways. The experimentally generated waves by walkers and superwalkers have been modeled using Bessel functions along with a Gaussian spatial decay~\citep{Molacek2013DropsTheory,molacek_bush_2013,tadrist_shim_gilet_schlagheck_2018,superwalkernumerical}. It would be interesting to explore this system in both $1$D and $2$D using a more realistic wave form that captures both oscillations and spatial decay of the wave field as observed in experiments. We still expect to observe the complex dynamical behaviors and interference effects, however, they may be attenuated in the parameter space due to the presence of spatial decay. It would also be interesting to explore the system presented in this work experimentally using walking and superwalking droplets. Although a similar setup has been explored by \citet{Braggdroplet2019}, the realization of the periodic potential using periodically spaced subsurface barriers alters the Faraday instability and hence the form of the droplet-generated waves. Conversely, the system studied in this work applies the periodic potential to the particle without affecting the droplet-generated waves. \citet{Perrard2014a} in their experiments of a walking droplet in a harmonic potential used magnetic fields to apply a harmonic potential to the particle alone without affecting the wave field. Perhaps, such a setup can also be used to generate a sinusoidal potential that acts only on the particle. Lastly, it would also be interesting to explore the system studied in this paper by adding disorder to the underlying periodic potential. It is known that quantum mechanical interference can result in localization of particles, known as Anderson localization, in a disordered potential~\citep{Andersonlocalization}. Such a localization effect may also emerge in our dynamical system upon adding disorder.      

\begin{acknowledgments}
J.P. was partly supported by Adelaide Summer Research Scholarship (ASRS) awarded by the University of Adelaide. R.V. was supported by Australian Research Council (ARC) Discovery Project DP200100834 during the course of the work. Some of the numerical results were computed using supercomputing resources provided by the Phoenix HPC service at the University of Adelaide.
\end{acknowledgments}

\section*{Data Availability Statement}

The data that support the findings of this study are available from the corresponding author upon reasonable request.

\appendix


\section{Dynamics of the low-memory system}\label{App: low memory}

The memory force term in Eq.~\eqref{Eq: dimless eq} is given by 
\begin{equation*}
F_M(t)=r\int_{-\infty}^t f(x_d(t)-x_d(s))\,\text{e}^{-(t-s)}\,\text{d}s.    
\end{equation*}
By making a change of variables $z=t-s$, we can transform this integral to
\begin{equation*}
F_M(t)=r\int_{0}^{\infty} f(x_d(t)-x_d(t-z))\,\text{e}^{-z}\,\text{d}z.    
\end{equation*}
Since the bouncing period of the droplet is $T_F$, the most recent wave was generated at time $t-T_F$ in dimensional variables which corresponds to time $t-1/\text{Me}$ in dimensionless variables (scaled by $T_F \text{Me}$). Hence, only considering the most recent impact, we can approximate the memory force as
\begin{equation*}
F_M(t)\approx r\,f(x_d(t)-x_d(t-1/\text{Me}))\,\text{e}^{-\frac{1}{\text{Me}}} \approx r\,f(\dot{x}_d(t)/\text{Me})\,\text{e}^{-\frac{1}{\text{Me}}}.    
\end{equation*}
Moreover, at low-memory, we have~\citep{Oza2013} $\text{Me} \sim {O}(1)$ so taking $\text{Me}=1$ we have the approximate form for the memory force
\begin{equation*}
F_M(t)\approx \frac{r}{e}\,f(\dot{x}_d(t)),
\end{equation*}
resulting in the equation of motion given in Eq.~\eqref{eq: lowmemode} and the system of nonlinear ODEs in Eq.~\eqref{eq: lowmem}.

\begin{figure}
\centering
\includegraphics[width=\columnwidth]{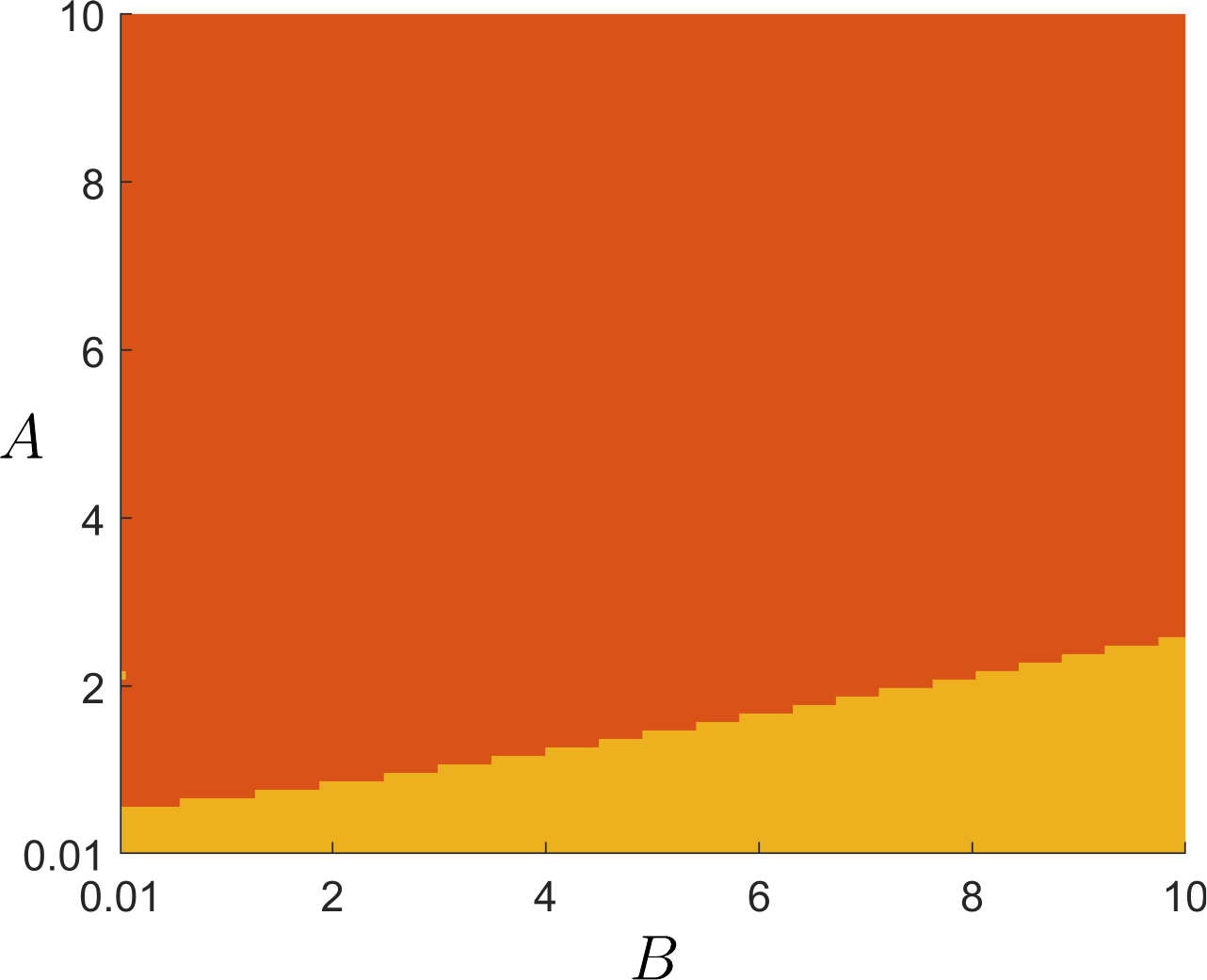}
\caption{Dynamical behaviors in the $(A, B)$ parameter space for the low-memory system given by Eq.~\eqref{eq: lowmem}. Here, red indicates back-and-forth oscillating walkers and yellow indicates runaway oscillating walkers. The parameters $\sigma=5$ and $r=5$ and the initial condition was fixed to $(x_d(0),v(0))=(\pi/B,0)$.}
\label{Fig: low memory}
\end{figure}

For the low-memory system  in Eq.~\eqref{eq: lowmem}, we find the same stationary states as the full system: $x_0=k\pi/B$ for $k\in \mathbb{Z}$. To determine the stability of the stationary solution in the
low-memory system, we perform a linear stability analysis by applying a small perturbation to the stationary solution in a similar way to that presented in Sec.~\ref{sec: lin stab full model}.

By applying a small perturbation to the equilibrium solution as follows: $(x_d, v)=({k\pi}/{B},0)+\epsilon(x_{d1},v_1)$, we obtain the linear system that governs the evolution of perturbations 

\begin{gather*}
     \begin{bmatrix} 
     \dot{x}_{d1} \\
     \dot{v}_1 \\
     \end{bmatrix}
     =
      \begin{bmatrix}
    0 & 1 \\
    (-1)^k AB\sigma & \sigma\left(\frac{r}{e}-1\right) \\
     \end{bmatrix}
      \begin{bmatrix}
      x_{d1}\\
      {v}_1
     \end{bmatrix}.
 \end{gather*}
 This yields the characteristic polynomial
 \begin{equation*}
     \lambda^2+\lambda\sigma\left(1-\frac{r}{e}\right)-(-1)^k\sigma AB=0.
 \end{equation*}
 Solving this quadratic gives us the eigenvalues
\begin{equation}\label{eq: lowmem_eigenvalues}
\lambda_{1,2}=\frac{\sigma}{2}\left( \frac{r}{e}-1 \right) \pm \sqrt{\left[ \frac{\sigma}{2}\left( \frac{r}{e}-1 \right) \right]^2 + (-1)^k A B \sigma}
 \end{equation}
From Eq.~\eqref{eq: lowmem_eigenvalues} we can see that for even $k$ we will have $\lambda_1>0$ and $\lambda_2<0$ and hence these stationary states are always unstable.

For odd $k$ we have two cases: (i) if $r<e$ then both the eigenvalues will have a negative real part making the stationary state stable, and (ii) if $r>e$ then both the eigenvalues will have a positive real part making the stationary state unstable. Thus the line $r=e$ separates the stable stationary states from the unstable stationary states.

To explore the unsteady dynamical behaviors, the system in Eq.~\eqref{eq: lowmem} is solved using the MATLAB ode45 solver. Fig.~\ref{Fig: low memory} shows the observed behaviors in a typical $(A,B)$ parameter-space plot for $\sigma=r=5$. Since $r=5>e$, the stationary state is unstable and we only see unsteady behaviors. In the unsteady regime, we find two types of behaviors for the WPE: (i) back-and-forth oscillating walkers (red) and (ii) runaway oscillating walkers (yellow). 

\bibliography{aipsamp}

\end{document}